\documentclass[10pt]{amsart}

\usepackage[T1]{fontenc}
\usepackage[utf8]{inputenc}
\usepackage[nomath]{lmodern}
\usepackage[expansion=false]{microtype}
\usepackage[
  paperwidth=170mm,
  paperheight=230mm,
  total={125mm,170mm},
  top=29mm,
  left=23mm,
  includeheadfoot
]{geometry}
\usepackage{xurl}
\usepackage{amsmath,amssymb,amsfonts,amsthm}
\usepackage{mathtools}
\usepackage{tikz}
\usepackage{tikz-cd}
\usetikzlibrary{arrows.meta,positioning,calc,fit,decorations.pathreplacing}
\usepackage{array}
\usepackage{enumitem}
\usepackage{xcolor}
\usepackage{hyperref}
\usepackage{orcidlink}
\usepackage{doi}

\hypersetup{
  hidelinks,
  pdfauthor={J. Rogelio Perez-Buendia, ORCID 0000-0002-7739-4779},
  pdftitle={An introduction to discrete dynamics in condensed sets: orbits, periodic loci, and arithmetic examples}
}

\newcommand{\ProFin}{\mathsf{ProFin}}
\newcommand{\Fin}{\mathsf{Fin}}
\newcommand{\Set}{\mathsf{Set}}
\newcommand{\Cond}{\mathsf{Cond}}
\newcommand{\Sh}{\mathsf{Sh}}
\newcommand{\Sub}{\operatorname{Sub}}
\newcommand{\Cov}{\operatorname{Cov}}
\newcommand{\Eq}{\operatorname{Eq}}
\newcommand{\Coeq}{\operatorname{Coeq}}
\newcommand{\Cont}{\operatorname{Cont}}
\newcommand{\Fix}{\operatorname{Fix}}
\newcommand{\Per}{\operatorname{Per}}
\newcommand{\PrePer}{\operatorname{PrePer}}
\newcommand{\Hist}{\operatorname{Hist}}
\newcommand{\admR}{\mathcal R}
\newcommand{\Z}{\mathbb{Z}}
\newcommand{\Q}{\mathbb{Q}}
\newcommand{\R}{\mathbb{R}}
\newcommand{\C}{\mathbb{C}}
\newcommand{\N}{\mathbb{N}}
\newcommand{\Ndisc}{\mathbb{N}_{\mathrm{disc}}}
\newcommand{\Zp}{\mathbb{Z}_p}
\newcommand{\Qp}{\mathbb{Q}_p}
\newcommand{\Qpalg}{\overline{\Q}_p}
\newcommand{\Cp}{\mathbb{C}_p}

\theoremstyle{plain}
\newtheorem{theorem}{Theorem}[section]
\newtheorem{proposition}[theorem]{Proposition}
\newtheorem{lemma}[theorem]{Lemma}

\theoremstyle{definition}
\newtheorem{definition}[theorem]{Definition}
\newtheorem{example}[theorem]{Example}
\newtheorem{remark}[theorem]{Remark}

\makeatletter
\renewenvironment{proof}[1][\proofname]{\par
  \pushQED{\qed}%
  \normalfont \topsep6\p@\@plus6\p@\relax
  \trivlist
  \item[\hskip\labelsep\itshape #1\@addpunct{:}]\ignorespaces
}{%
  \popQED\endtrivlist\@endpefalse
}

\def\ps@headings{\ps@empty
  \def\@evenhead{%
    \setTrue{runhead}%
    \normalfont\scriptsize
    \rlap{\thepage.}\hfil
    \def\thanks{\protect\thanks@warning}%
    \leftmark{}{}\hfil}%
  \def\@oddhead{%
    \setTrue{runhead}%
    \normalfont\scriptsize \hfil
    \def\thanks{\protect\thanks@warning}%
    \rightmark{}{}\hfil \llap{\thepage.}}%
  \let\@mkboth\markboth
}
\makeatother
\begin{document}

\title[An introduction to discrete dynamics in condensed sets]{An introduction to discrete dynamics in condensed sets:\\
orbits, periodic loci, and arithmetic examples}

\author[J. Rogelio P\'erez-Buend\'ia]{J.~Rogelio P\'erez-Buend\'ia\,\orcidlink{0000-0002-7739-4779}}
\thanks{Published in \emph{Memorias de la Sociedad Matem\'atica Mexicana} \textbf{10} (2026), 87--122.}
\address{SECIHTI -- Centro de Investigaci\'on en Matem\'aticas (CIMAT),\newline
Unidad M\'erida, M\'erida, Mexico}
\email{rogelio.perez@cimat.mx}



\begin{abstract}
We formulate discrete dynamics internally in the topos $\Cond$ of condensed
sets. For an endomorphism $f:X\to X$, we construct orbit quotients, backward
histories, preimage trees, and periodic loci as objects evaluated on profinite
parameter spaces. The framework recovers global-section dynamics, fully
faithfully embeds compact Hausdorff systems, characterizes the invertible
extension, and is compatible with inverse limits. It organizes $p$-adic
congruence towers and shows, for a Cantor-type non-Archimedean Julia set, that
classical periodic points can be dense although the internal periodic locus is
not isomorphic to the compact condensed Julia set.
\end{abstract}

\maketitle

\begin{center}
\begin{minipage}{\textwidth}
\small
\noindent\textsc{Resumen.}
Formulamos la din\'amica discreta internamente en el topos $\Cond$ de conjuntos
condensados. Para un endomorfismo $f:X\to X$, construimos cocientes orbitales,
historias hacia atr\'as, \'arboles de preim\'agenes y loci peri\'odicos como
objetos evaluados en espacios profinitos de par\'ametros. El formalismo recupera
la din\'amica en secciones globales, encaja plenamente y fielmente los sistemas
compactos Hausdorff, caracteriza la extensi\'on invertible y es compatible con
l\'imites inversos. Organiza torres de congruencias $p$-\'adicas y muestra, para
un conjunto de Julia no arquimediano tipo Cantor, que los puntos peri\'odicos
cl\'asicos pueden ser densos aunque el locus peri\'odico interno no sea isomorfo
al conjunto condensado compacto de Julia.
\end{minipage}
\end{center}

\section{Introduction}
\label{sec:intro}

\begin{center}
\begin{minipage}{0.72\textwidth}
\small\itshape
``L'autre est celle de la mer... Pourtant elle finit par entourer la
substance r\'etive.''

\hfill A.~Grothendieck, \emph{R\'ecoltes et Semailles}
\end{minipage}
\end{center}

Finite dynamics is especially concrete. A self-map can be represented by its
functional graph, whose cycles display periodic behavior directly. See, for
instance, \cite{MitchellMukherjeeVojtechovskyEndofunctions}. Two phenomena that
arise beyond a single finite level motivate the constructions developed below.

Arithmetic and non-Archimedean dynamics enrich this picture with compatible
structure across scales. A map over the \(p\)-adic integers \(\Zp\), with
\(p\) prime, is seen through all of its reductions modulo \(p^m\), so its
periodic data form a system across a congruence tower. Already the translation
\(x\mapsto x+1\) shows the point: on the finite quotient
\(\Z/p^m\Z\), every point is periodic of exact period \(p^m\), whereas the
inverse-limit system on \(\Zp\) has no periodic points.

In non-Archimedean geometry, periodic points may also be dense in a compact
Cantor-type Julia set, while the dense periodic part and the compact ambient
space carry different structures. These examples point to the same need:
periodic data should be kept together with their variation in families, their
finite-level compatibility and their position inside compact or profinite phase
spaces.

Condensed mathematics provides a language for these compatibilities.
Introduced by Clausen--Scholze, it organizes objects whose algebraic structure
and topology interact as sheaves on profinite sets
\cite{ClausenScholze,LeStumCondensed}. For the purposes of this article, a
condensed set \(X\) assigns to each profinite set \(S\) a set \(X(S)\) of
\(S\)-parameterized families, subject to gluing along finite jointly surjective
covers. The resulting category \(\Cond\) is a Grothendieck topos. Thus it is
equivalent to a category of sheaves on a site, that is, on a category equipped
with a specified notion of covering. Evaluation at the one-point profinite set
gives the ordinary points, or global sections, \(X(*)\), while a general \(S\)
supplies the family data used below.

To the author's knowledge, this is the first systematic formulation of
condensed discrete dynamics in the precise sense used here: the study of
endomorphisms and their orbit and return constructions internally in the topos
\(\Cond\). A condensed dynamical system is an endomorphism \(f:X\to X\) in
\(\Cond\), equivalently an action of the discrete condensed monoid associated
with \(\mathbb N\).
From this datum we construct orbit quotients, backward histories, preimage
trees, and fixed, periodic and preperiodic loci.

The resulting language recovers ordinary dynamics on global sections and
embeds compact Hausdorff dynamical systems fully faithfully. Backward histories
satisfy the universal property of the invertible extension, while periodic loci
commute with inverse limits in \(p\)-adic congruence towers. The geometric
comparison appears in a Cantor-type non-Archimedean Julia example over the
complete algebraically closed \(p\)-adic field \(\Cp\). The internal periodic
locus is the discrete condensed object associated with the periodic points. Its
global points are dense in the classical Julia set, and the internal locus is
not isomorphic in \(\Cond\) to the compact condensed Julia set.

The paper is organized as follows. The remainder of the Introduction fixes the
notation and locates the internal layer among related categorical and observed
frameworks. Section~\ref{sec:condensed} recalls the profinite site and condensed
sets. Sections~\ref{sec:systems} and~\ref{sec:periodic} construct the orbit,
history and return objects, including their behavior in \(p\)-adic congruence
towers. Section~\ref{sec:not-counting} compares discrete and compact periodic
structures and applies the comparison to a Cantor-type Julia set.

\subsection{Notation and conventions}
\label{sec:notation}

We use the following notation. Throughout, $p$ denotes a prime
number, and
\[
  \Zp=\varprojlim_m \Z/p^m\Z
\]
is the ring of $p$-adic integers. Equivalently, $\Zp$ is obtained by completing
$\Z$ for the metric induced by the $p$-adic absolute value. Appendix~\ref{app:topological-background}
recalls this metric viewpoint. Its fraction field is $\Qp$, the field of
$p$-adic numbers, and $\Cp$ is the completion of an algebraic closure of $\Qp$.
The notation $\Ndisc$ denotes the condensed set associated with the discrete
space $\N$. The symbol $\ProFin$ denotes the category of profinite sets and
continuous maps. It is introduced carefully in Section~\ref{sec:condensed}. We
write $\Set$ for the category of sets, $\Cont(S,T)$ for continuous maps, and
$\Eq$ and $\Coeq$ for equalizers and coequalizers. When a topological space is
converted into a condensed set, we work in the standard convenient category of
compactly generated weak Hausdorff spaces \cite{ClausenScholze,LeStumCondensed}.
For the examples below, the reader may safely keep compact metrizable
spaces, profinite spaces, and finite discrete spaces in mind.

\subsection{Relation with categorical and observed dynamics}
\label{sec:relation}

Before constructing the dynamical objects, we use \emph{internal} and
\emph{external} in the following methodological sense. An internal construction
is built from the categorical operations available in the
ambient topos: limits, colimits, pullbacks, subobjects, joins and inverse images.
An \emph{external} or \emph{observed} layer enriches the internal
datum by extra choices not determined by the ambient category alone. Those
choices may still be made out of objects already present in the category, such
as subobjects or covering families. They enter as additional structure, just as
a topology selects a family of subsets on an underlying set.

More generally, in any category $\mathcal C$, iteration is encoded by an
endomorphism $f:X\to X$. This is the categorical syntax used below: conjugacy,
factors, subsystems and invariance can all be expressed in terms of morphisms.
For surrounding categorical approaches to dynamics and systems, see
\cite{LawvereToposes,MyersDouble}. Sheaf-theoretic and
topos-theoretic treatments of discrete dynamics include
\cite{SchultzSpivakVasilakopoulou,GarofaliDynamicalSheaves,HoraKamio,RogersTopologicalMonoid}.
In this article, that categorical syntax is interpreted in $\Cond$ through
evaluation on profinite sets.

Appendix~\ref{app:categorical-background} gives the precise definitions of a
site and a Grothendieck topos and recalls the categorical properties needed
later.
In particular, the subobjects of any object in a Grothendieck topos $\mathcal E$
form a complete lattice \cite{MacLaneMoerdijk,JohnstoneTopos}. For an object
$X$, we write $\Sub(X)$ for this lattice of subobjects when the ambient topos is
clear. In $\Cond$, we write $\Sub_{\Cond}(X)$ when the ambient category
matters.

A different, complementary construction is the observed Fatou--Julia framework
developed in \cite{PerezBuendiaFatuTopos}. This is
external to the internal layer in the sense just described: one enriches the
dynamical object by choosing an observation regime,
\[
  (X,f,\mathcal O,\admR),
\]
where $f:X\to X$ is the dynamical endomorphism in that observed layer,
$\mathcal O\subseteq \Sub(X)$ is an observable subframe and
$\admR\subseteq\Cov(\mathcal O)$ is a class of admissible covers, with
$\Cov(\mathcal O)$ denoting covering families in the observable frame. The extra
data are used to formulate cover-wise equicontinuity, observable Fatou loci,
observable Julia loci, and coreflective universal properties inside the chosen
observable regime. This is the layer in which point-free Fatou--Julia questions
are naturally formulated.

\section{Profinite sets, sites and sheaves}
\label{sec:condensed}

\subsection{Profinite sets as finite-layer spaces}

Profinite sets encode compatible finite-level data and provide the parameter
spaces on which condensed sets are evaluated.

Let $\Fin$ denote the category of finite discrete sets.
A \emph{profinite set} is a compact Hausdorff totally disconnected space.
Equivalently, it is a cofiltered inverse limit of finite discrete sets:
\[
  S\cong\varprojlim_i S_i,
  \qquad S_i\in\Fin.
\]
We write $\ProFin$ for the category of profinite sets and continuous maps.
The equivalence between the topological and inverse-limit descriptions makes
profinite sets useful in arithmetic dynamics: a point of the compact space $S$
is read through its compatible images in the finite quotients $S_i$
\cite{JohnstoneStoneSpaces}. This finite-level description is the bridge between
Stone duality, $p$-adic arithmetic, finite quotient dynamics and the condensed
site.
For an elementary
introduction from an arithmetic viewpoint, see the published expository account
\cite{PerezBuendiaStoneGaceta}. A longer self-contained version is
\cite{PerezBuendiaStoneArxiv}.

The central arithmetic example is
\[
  \Zp=\varprojlim_m\Z/p^m\Z,
\]
where $p$ is a fixed prime. An element of $\Zp$ is a compatible system of
residue classes $(x_m)_m$ with
$x_{m+1}\equiv x_m\pmod {p^m}$.
A $1$-Lipschitz endomorphism of $\Zp$, for instance a polynomial map
with coefficients in $\Zp$, induces compatible maps on the finite quotients
$\Z/p^m\Z$ and therefore a compatible system of finite dynamical systems.

There is an equivalent analytic construction. Appendix~\ref{app:topological-background}
recalls the $p$-adic metric and the topology of $\Qp$ explicitly. Completing
$\Z$ for that metric gives the same profinite ring $\Zp$. Completing $\Q$ for
the associated non-Archimedean absolute value gives $\Qp$, and the algebraic
description agrees through the canonical identification
$\Qp\cong\operatorname{Frac}(\Zp)$. Passing to an algebraic closure and
completing gives the usual field $\Cp$ for one-variable non-Archimedean
dynamics. Standard background can be found in
\cite{GouveaPadic,SilvermanADS,BenedettoNAM}. Here $\Qpalg$ denotes a chosen
algebraic closure of $\Qp$. The
associated chain is
\[
  \Zp\subset \Qp\subset \Qpalg\subset \Cp
\]
and the absolute value extends along it.

\begin{figure}[htbp]
\centering
\begin{tikzcd}[row sep=2.7em, column sep=3.0em]
& \Zp
    \arrow[dl, "\pi_1"']
    \arrow[d,  "\pi_2" description]
    \arrow[dr, "\pi_3"]
    \arrow[drr, dashed, bend left=12, "\pi_m"]
  & & \\
\Z/p\Z
  & \Z/p^2\Z \arrow[l, "\bmod p"']
  & \Z/p^3\Z \arrow[l, "\bmod p^2"']
  & \cdots \arrow[l]
\end{tikzcd}
\caption{$\Zp$ as the inverse limit of its finite congruence levels.}
\label{fig:zp-tower}
\end{figure}
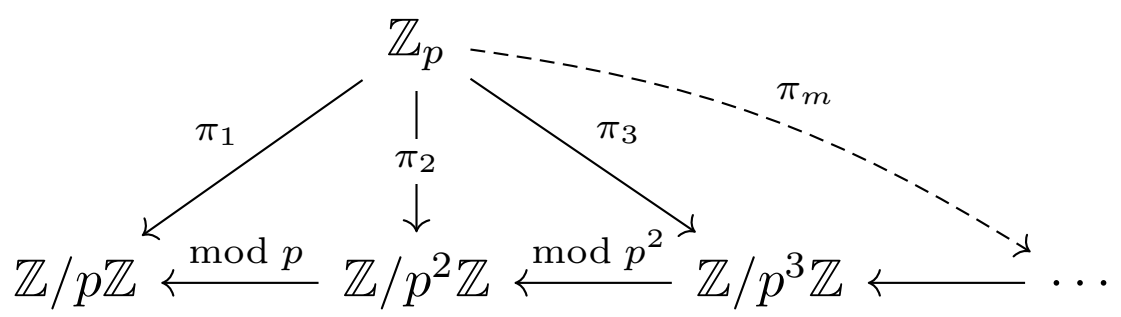

\subsection{Sites and sheaves: the gluing idea}

A \emph{site}, for our purposes, is a category together with a specified class
of covering families. A \emph{sheaf} on a site is a contravariant functor that
glues compatible local data along those coverings. This is the level of site
theory needed in the main text. Standard references are
\cite{MacLaneMoerdijk,JohnstoneTopos}. Appendix~\ref{app:categorical-background}
summarizes the categorical conventions used below.

The category $\ProFin$ is equipped with the finite jointly surjective
pretopology used in the profinite-site presentation of condensed sets
\cite{ClausenScholze,LeStumCondensed}: a finite family $\{U_i\to S\}$ covers
$S$ when $\coprod_i U_i\to S$ is surjective. These profinite sets are the
objects on which condensed sets are evaluated.

The sheaf condition has the same geometric meaning here as in ordinary
topology. If $S$ is covered by finitely many profinite pieces $U_i$, then a
section over $S$ is the same as a family of sections over the $U_i$ that agree
on all overlaps $U_i\times_S U_j$. The fiber product $U_i\times_S U_j$ plays
the role of the overlap $U_i\cap U_j$ in an ordinary open cover. The difference
is not the gluing idea but the site: instead of open subsets of one fixed
space, the covers live among profinite sets themselves.

For a sheaf $X$ on this site, the formal expression of this gluing rule is the
equalizer diagram
\[
  X(S)\longrightarrow \prod_i X(U_i)
  \rightrightarrows \prod_{i,j}X(U_i\times_S U_j)
\]
for every finite jointly surjective cover $\{U_i\to S\}$. The first arrow
restricts a global section over $S$ to the pieces $U_i$. The two arrows on the
right compare the two induced restrictions to each overlap $U_i\times_S U_j$.
Saying that the diagram is an equalizer means that $X(S)$ is exactly the set of
compatible families on the pieces.

\begin{figure}[htbp]
\centering
\begin{tikzpicture}[
  node distance=1.3cm and 1.6cm,
  obj/.style={draw, rounded corners=2pt, inner sep=4pt, font=\small},
  arr/.style={-{Latex[length=2mm]}, thick}
]
\node[obj] (S) {$S$};
\node[obj, above left=of S] (U1) {$U_1$};
\node[obj, above=of S] (U2) {$U_2$};
\node[obj, above right=of S] (Uk) {$U_r$};
\node[obj, below=of S] (coprod) {$\coprod_i U_i$};
\draw[arr] (U1) -- (S);
\draw[arr] (U2) -- (S);
\draw[arr] (Uk) -- (S);
\draw[arr] (coprod) -- node[right, font=\scriptsize] {surjective} (S);
\node[font=\scriptsize, above=0.15cm of U2] {finite family in $\ProFin$};
\node[font=\scriptsize, below=0.15cm of coprod] {covering condition};
\end{tikzpicture}
\caption{A finite jointly surjective cover in the condensed site.}
\label{fig:profin-cover}
\end{figure}
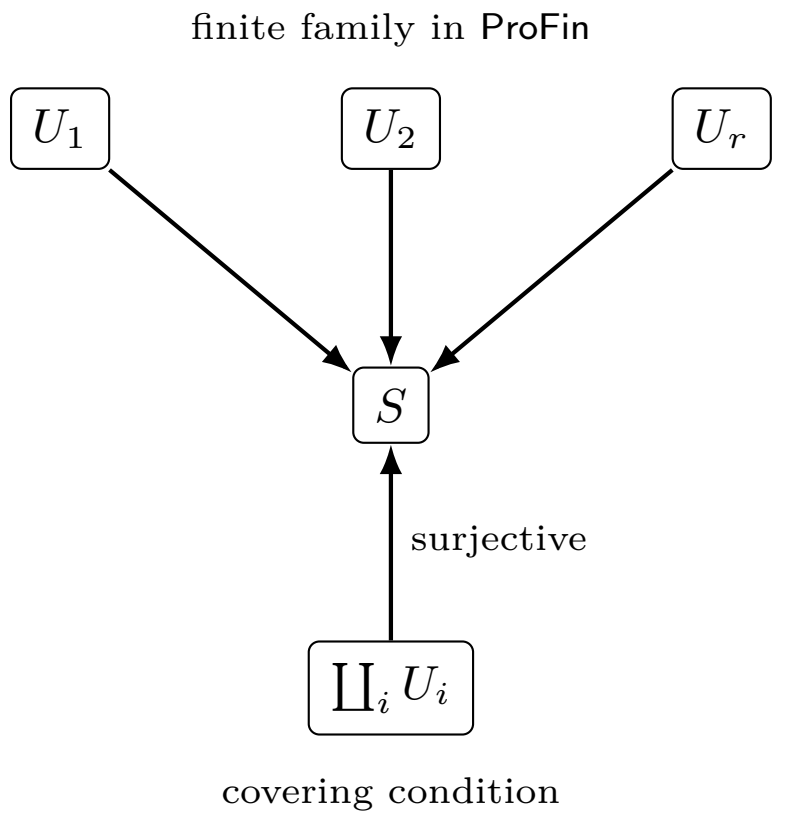

\begin{figure}[htbp]
\centering
\begin{tikzcd}[column sep=large]
X(S) \arrow[r] &
\prod_i X(U_i)
\arrow[r, shift left=0.55ex]
\arrow[r, shift right=0.55ex, swap] &
\prod_{i,j}X(U_i\times_S U_j)
\end{tikzcd}
\[
\begin{array}{c}
\text{global section on }S\\
\Updownarrow\\
\text{compatible sections on the finite cover }\{U_i\to S\}
\end{array}
\]
\caption{The sheaf condition for a finite profinite cover.}
\label{fig:sheaf-equalizer}
\end{figure}

\subsection{Condensed sets}

\begin{definition}[Condensed set]
Fix a set-theoretic universe, or equivalently a small skeleton of the profinite
site. In the working presentation used here, a \emph{condensed set} is a functor
\[
  X:\ProFin^{\mathrm{op}}\to\Set
\]
satisfying the above gluing condition for every finite jointly surjective
covering in $\ProFin$. Let $\tau_{\mathrm{fs}}$ denote this finite jointly
surjective pretopology. We write
\[
  \Cond=\Sh(\ProFin,\tau_{\mathrm{fs}})
\]
for the resulting category of condensed sets. This category is a Grothendieck
topos \cite{ClausenScholze,LeStumCondensed}.
\end{definition}

We use this presentation throughout. Foundational treatments
often use equivalent small sites, for instance sites built from extremally
disconnected profinite sets with the coherent topology. The internal dynamical
constructions below do not depend on that choice
\cite{ClausenScholze,LeStumCondensed}.

\[
  (\ProFin,\tau_{\mathrm{fs}})\xrightarrow{\ \Sh\ }\Cond,
  \qquad
  X\in\Cond\longmapsto X(S)\in\Set .
\]

The set $X(S)$ consists of the $S$-parameterized sections of $X$.
When $S$ changes, functoriality gives restriction maps along continuous maps of
profinite sets.

If $T$ is a compactly generated weak Hausdorff space, its associated condensed
set is
\[
  X_T(S)=\Cont(S,T).
\]
At the terminal profinite set $*$ one has $X_T(*)\cong T$. At a nontrivial
profinite set $S$, one obtains continuous $S$-families of points of $T$.
For compact Hausdorff spaces, and in particular for profinite spaces, this
construction faithfully remembers the topology
\cite{ClausenScholze,LeStumCondensed}.
When no ambiguity can arise, we use the same symbol for a profinite space and
for its associated condensed set.
The compactly generated weak Hausdorff hypothesis is a standard technical
choice that gives a well-behaved category of spaces for mapping constructions.
All compact Hausdorff, finite discrete, profinite, and compact metrizable
spaces considered here belong to this class.

\begin{example}[Finite and profinite targets]
If $F$ is a finite discrete set, then $X_F(S)=\Cont(S,F)$ is the set of locally
constant maps $S\to F$. Equivalently, each such map factors through some finite
quotient $S\twoheadrightarrow S_i$:
\[
  \Cont(S,F)\cong
  \varinjlim_{(S\twoheadrightarrow S_i)}\operatorname{Map}(S_i,F).
\]
Thus a finite discrete set is remembered by all locally constant labelings of
finite clopen partitions of $S$.

If $P$ is itself profinite, then $X_P(S)=\Cont(S,P)$ contains genuinely
continuous profinite families. For instance, $X_{\Zp}(\Zp)$ contains the
identity map $\Zp\to\Zp$, which does not factor through any finite quotient of
the source: if it factored through a finite quotient, its image would be finite,
contrary to the infinitude of $\Zp$.
\end{example}

\begin{remark}
The terminology follows the usual sheaf convention. If $\mathcal F$ is a sheaf
on a topological space $T$, the category $\operatorname{Open}(T)$ has terminal
object $T$, and the global sections are $\mathcal F(T)$. The profinite site has
terminal object $*$, the one-point profinite set. Its representable sheaf is
terminal in $\Cond$, so the Yoneda lemma gives
\[
  \Gamma(X)=\operatorname{Hom}_{\Cond}(1,X)\cong X(*).
\]
For $X=X_T$ associated with a topological space $T$, this recovers
$X_T(*)\cong T$. Evaluation on a general profinite set $S$ gives
$S$-parameterized families. Terminal evaluation is therefore evaluation at the
terminal object, rather than a stalk construction.
\end{remark}

The next lemma passes from filtered joins of subobjects to finite-stage local
data. It will be used for periodic loci and discrete condensed periodic objects.

\begin{lemma}[Finite-stage detection for filtered unions]
\label{lem:finite-stage-detection}
Let $\{A_i\}_{i\in I}$ be a filtered system of subobjects of a condensed set
$X$, with transition maps monomorphisms, and let
\[
  A=\bigvee_{i\in I}A_i\hookrightarrow X
\]
be their join in $\Sub_{\Cond}(X)$. If $S$ is profinite and $s\in A(S)$, then
there is an index $i\in I$ and a section $t\in A_i(S)$ whose image in $X(S)$ is
$s$.
\end{lemma}

\begin{proof}
The join is the sheafification of the filtered presheaf union inside $X$.
Thus a section of $A$ over $S$ is locally represented by a section of the
presheaf colimit. The topology on $\ProFin$ uses finite jointly surjective
covers, so this local representation may be taken over a finite cover. On each
member of the cover, an element of a filtered colimit of sets is represented at
some finite stage. Filteredness gives one index dominating the finitely many
indices that occur. After replacing the local representatives by their images
in this common stage $A_i$, their restrictions to overlaps agree, because they
have the same image in $X$ and $A_i\to X$ is a monomorphism. Since $A_i$ is a
sheaf, these local sections glue to a section of $A_i(S)$, and its image in
$X(S)$ is $s$. This is the standard construction of colimits in a sheaf topos
\cite{MacLaneMoerdijk,JohnstoneTopos}.
\end{proof}

\section{Internal condensed dynamical systems and orbit constructions}
\label{sec:systems}

We now construct the orbit data attached to an endomorphism in $\Cond$ and
describe their evaluation on profinite parameter spaces.

\begin{definition}
A \emph{condensed dynamical system} is a pair $(X,f)$ where $X\in\Cond$ and
$f:X\to X$ is a morphism in $\Cond$.
\end{definition}

Equivalently, a condensed dynamical system is a monoid action
$\Ndisc\times X\to X$ in $\Cond$, where $\Ndisc$ is the discrete condensed
monoid associated with $\N$. Since $\N$ is freely generated by one element,
such an action is the same as an endomorphism $f:X\to X$. In
categorical terms it is a functor
\[
  B\N\longrightarrow \Cond,
\]
where $B\N$ is the one-object category whose endomorphism monoid is $\N$. The
unique object is sent to $X$, and the generator $1\in\N$ is sent to the
endomorphism $f:X\to X$. This realizes the action generated by $f$, with the
periodic loci as one part of that action.
The same construction extends to actions of more general monoid objects in
$\Cond$. We focus on the discrete-time case generated by one endomorphism.

For each $S\in\ProFin$, naturality gives a map
\[
  f_S:X(S)\to X(S).
\]
Thus a condensed dynamical system is a compatible family of dynamical systems
on the sets of $S$-sections. If $X=X_T$ comes from a topological dynamical
system $f:T\to T$, then $f_S$ sends a continuous map $\alpha:S\to T$ to
$f\circ\alpha$.

\begin{figure}[htbp]
\centering
\[
\begin{tikzcd}[column sep=large, row sep=large]
X(S) \arrow[r, "f_S"] \arrow[d, "\theta^\ast"'] &
X(S) \arrow[d, "\theta^\ast"]\\
X(S') \arrow[r, "f_{S'}"'] &
X(S')
\end{tikzcd}
\qquad
\begin{tikzcd}[column sep=large]
S' \arrow[r, "\theta"] & S
\end{tikzcd}
\]
\caption{Naturality of evaluation on profinite sets.}
\label{fig:evaluation-natural}
\end{figure}

The forward action morphism associated with $f$ is
\[
  \operatorname{orb}^+_f:\Ndisc\times X\longrightarrow X,
  \qquad
  (n,x)\longmapsto f^n(x).
\]
Evaluating on a profinite set $S$ gives
\[
  \operatorname{orb}^+_{f,S}:
  \Ndisc(S)\times X(S)\longrightarrow X(S).
\]
An element of
$\Ndisc(S)=\Cont(S,\N)$ is a locally constant choice of time on $S$. Its image is
finite by compactness of $S$, so it is constant on the pieces of a finite clopen
partition. Thus the condensed action morphism allows the iterate $n$ to vary
locally constantly over the profinite parameter space. For constant time
functions, this reduces to the ordinary iterate $f_S^n$.

Classically, one forms the quotient of a set by the orbit relation. In the
condensed setting the quotient is formed internally, so the resulting object is
still a sheaf of compatible orbit families. The forward orbit object is the
internal
coequalizer
\[
\begin{tikzcd}
X \arrow[r, shift left=0.45ex, "\mathrm{id}_X"]
  \arrow[r, shift right=0.45ex, swap, "f"] &
X \arrow[r] &
X_{\N}.
\end{tikzcd}
\]
Equivalently,
\[
  X_{\N}:=\Coeq(\mathrm{id}_X,f).
\]
The subscript $\N$ indicates quotient by the forward $\N$-action generated by
$f$.
In sets this is the quotient by the equivalence relation generated by
$x\sim f(x)$.

\begin{remark}[Orbit quotients are sheaf quotients]
In $\Cond$ the object $X_{\N}$ is the internal quotient in the topos. Therefore
one should not identify $X_{\N}(S)$ naively with the pointwise quotient
$X(S)/{\sim}$ without hypotheses: colimits of sheaves include sheafification.
See Appendix~\ref{app:categorical-background} and
\cite{MacLaneMoerdijk,JohnstoneTopos}. This is one of the places where the
condensed setting is structural. The orbit quotient is a quotient of compatible
families formed with the sheaf condition still present.
\end{remark}

For backward histories and natural extensions, the noninvertible case is the
one that matters. A backward orbit is not usually determined by one point. It
is a compatible choice of preimages. Internally this is the inverse limit
\[
  \Hist^{-}(X,f)
  =
  \varprojlim
  \left(
  \cdots \xrightarrow{f} X \xrightarrow{f} X
  \right).
\]
Since limits of sheaves are computed sectionwise
\cite{MacLaneMoerdijk,JohnstoneTopos}, evaluation on $S$ gives
\[
  \Hist^{-}(X,f)(S)
  =
  \left\{
  (x_i)_{i\ge0}\in X(S)^{\N}:
  f_S(x_{i+1})=x_i\ \text{for all }i\ge0
  \right\}.
\]
The $S$-sections of $\Hist^{-}(X,f)$ are therefore compatible histories
parameterized by $S$. If $f$ is an automorphism, one may instead use the
associated $\Z$-action to describe bilateral orbits.

\begin{proposition}[Universal natural extension]
\label{prop:history-universal}
Let $(X,f)$ be a condensed dynamical system. The object
$\Hist^{-}(X,f)$ carries a canonical automorphism $\sigma$, and the projection
$\pi_0:\Hist^{-}(X,f)\to X$ satisfies
\[
  \pi_0\circ \sigma=f\circ \pi_0.
\]
The triple $(\Hist^{-}(X,f),\sigma,\pi_0)$ is terminal in the following sense:
for every triple $(Y,g,u)$ where $g:Y\to Y$ is an automorphism in $\Cond$ and
$u:Y\to X$ satisfies $u\circ g=f\circ u$, there is a unique morphism
$\widetilde u:Y\to\Hist^{-}(X,f)$ with
\[
  \pi_0\circ\widetilde u=u,
  \qquad
  \widetilde u\circ g=\sigma\circ\widetilde u.
\]
\end{proposition}

\begin{proof}
A section of $\Hist^{-}(X,f)$ over a profinite set $S$ is a sequence
$(x_0,x_1,x_2,\ldots)$ with $f_S(x_{i+1})=x_i$ for all $i\ge0$. Define
\[
  \sigma_S(x_0,x_1,x_2,\ldots)=(f_S(x_0),x_0,x_1,\ldots).
\]
Its inverse is
\[
  \sigma_S^{-1}(x_0,x_1,x_2,\ldots)=(x_1,x_2,\ldots).
\]
These maps are natural in $S$, hence define an automorphism $\sigma$ of
$\Hist^{-}(X,f)$.
For $k\ge0$, write $\pi_k:\Hist^{-}(X,f)\to X$ for the $k$-th projection.

Let $(Y,g)$ be a condensed dynamical system with $g$ an automorphism, and let
$u:Y\to X$ satisfy $u\circ g=f\circ u$. For each $k\ge0$ set
\[
  u_k=u\circ g^{-k}:Y\to X.
\]
Then
\[
  f\circ u_{k+1}
  =f\circ u\circ g^{-(k+1)}
  =u\circ g\circ g^{-(k+1)}
  =u\circ g^{-k}
  =u_k.
\]
Thus the morphisms $(u_k)_{k\ge0}$ form a compatible cone over the inverse
system defining $\Hist^{-}(X,f)$, and therefore induce a unique morphism
\[
  \widetilde u:Y\to\Hist^{-}(X,f)
\]
with $\pi_k\circ \widetilde u=u_k$. The identity
$\widetilde u\circ g=\sigma\circ\widetilde u$ follows after composing with all
projections $\pi_k$, and uniqueness follows because the projections jointly
determine morphisms into a limit.
\end{proof}

For preimage trees, let $a:S\to X$ be an $S$-parameterized point of $X$, where
the profinite set $S$ is identified with its representable condensed set. The
$n$-th level of the preimage tree over $a$ is the pullback in the slice over $S$
\[
  T_n(a):=S\times_{X,\,a,\,f^n}X.
\]
Equivalently, a section over a further profinite set $T\to S$ is a section
$x\in X(T)$ satisfying
\[
  f_T^n(x)=a|_T.
\]
The morphism $f$ induces transition maps $T_{n+1}(a)\to T_n(a)$, so these
levels form a tower over $S$.
These pullback towers are the internal version of the preimage trees familiar
in arithmetic dynamics \cite{SilvermanADS}. When $X$ comes from a variety or
analytic space over a field and $a$ is a rational point, the finite levels
$f^{-n}(a)$ often carry natural Galois actions. Organizing those actions leads
to orbital arboreal representations and to dynamical
N\'eron--Ogg--Shafarevich-type criteria. For arithmetic instances of this
viewpoint, see \cite{PerezBuendiaNOS}. Here this arithmetic motivation is used
to isolate the underlying internal preimage construction in $\Cond$, where the
whole compatible tower remains available for later constructions.

\begin{center}
\small
\begin{tabular}{>{\raggedright\arraybackslash}p{0.24\textwidth}|>{\raggedright\arraybackslash}p{0.33\textwidth}|>{\raggedright\arraybackslash}p{0.31\textwidth}}
Construction & Internal construction in $\Cond$ & Evaluation on $S$\\
\hline
Forward action morphism &
$\Ndisc\times X\to X$ &
$(n,x)\mapsto f_S^n(x)$, with $n$ locally constant on $S$\\
Orbit quotient &
$\Coeq(\mathrm{id}_X,f)$ &
Condensed quotient of orbit families\\
Backward history &
$\varprojlim(\cdots\xrightarrow f X\xrightarrow f X)$ &
Sequences $(x_i)$ with $f_S(x_{i+1})=x_i$\\
Preimage level over $a:S\to X$ &
$S\times_{X,\,a,\,f^n}X$ &
Families solving $f^n(x)=a$ after pullback
\end{tabular}
\end{center}

The first two columns use categorical constructions available in any complete
and cocomplete topos. The third gives their evaluation on profinite sets in
$\Cond$, including descent along finite jointly surjective covers.

\begin{definition}
Let $(X,f)$ and $(Y,g)$ be condensed dynamical systems. A morphism
$u:(X,f)\to(Y,g)$ is a morphism $u:X\to Y$ in $\Cond$ such that
$u\circ f=g\circ u$. An isomorphism $u$ gives a conjugacy. If $u$ is an
epimorphism, $(Y,g)$ is a factor of $(X,f)$. If $u$ is a monomorphism, $(X,f)$
is a subsystem of $(Y,g)$.
\end{definition}

\begin{definition}[Invariant subobjects]
\label{def:invariance}
Let $A\hookrightarrow X$ be a subobject of a condensed dynamical system
$(X,f)$.
\begin{enumerate}[label=(\alph*)]
\item $A$ is \emph{forward invariant} if the composite
$A\hookrightarrow X\xrightarrow{f}X$ factors through $A$.
\item $A$ is \emph{backward invariant} if the pullback
\[
  f^{-1}(A):=X\times_{X,f}A
\]
factors through $A$ as a subobject of $X$.
\item $A$ is \emph{invariant} if it is both forward and backward invariant.
\end{enumerate}
\end{definition}

For a representable profinite phase space $P$, this recovers the usual
conditions for a continuous map $f:P\to P$: forward invariance is the condition
$f(A)\subseteq A$, and backward invariance is the condition
$f^{-1}(A)\subseteq A$.

\begin{example}[Finite dynamics]
Let $F$ be a finite set, viewed as a finite discrete profinite set. Any map
$f:F\to F$ defines a condensed dynamical system. Evaluating at a finite profinite set
$S$ gives the product system on $\operatorname{Map}(S,F)$: an $S$-indexed
configuration evolves componentwise.
\end{example}

\begin{example}[The $p$-adic odometer]
Let $\tau:\Zp\to\Zp$ be translation by one, $\tau(x)=x+1$. On the finite
quotient $\Z/p^m\Z$, the induced map is a single cycle of length $p^m$. Thus
every finite quotient has only periodic points, but the required return time
grows with the level. In the inverse limit, $\tau$ has no point of positive
period: $\tau^n(x)=x$ would imply $n=0$ in the characteristic-zero ring $\Zp$.
Equivalently, for every fixed $n\ge 1$ the equalizer of $\tau^n$ and the identity
on $\Zp$ is empty, although the finite quotients have returns of unbounded
period. This example illustrates the information retained by fixed-iterate
equalizers across the congruence tower.
\end{example}

\begin{example}[Profinite shifts]
If $A$ is a finite alphabet, the one-sided shift on $A^{\N}$ with the product
topology is a continuous endomorphism of a profinite set. It therefore defines
a condensed dynamical system. The condensed viewpoint keeps track of orbits in
$A^\N$ together with continuous families of such orbits over every profinite set.
\end{example}

\begin{example}[Classical systems as condensed systems]
Let $T$ be a compactly generated weak Hausdorff space and let $f:T\to T$ be
continuous. The associated condensed system is
\[
  (X_T,X_f),\qquad X_T(S)=\Cont(S,T),
\]
where
\[
  X_f(S):\Cont(S,T)\to\Cont(S,T),\qquad
  \alpha\mapsto f\circ\alpha.
\]
Thus every ordinary topological dynamical system has a condensed image. At the
terminal profinite set one has $X_T(*)\cong T$, and the induced map on global
sections is $f$ itself. Thus the original topological dynamical system is
recovered on global sections. At a general $S$ one sees continuous
$S$-families of orbits.
\end{example}

For the full-faithfulness statement below, we restrict to compact Hausdorff
phase spaces.

\begin{proposition}[Full faithfulness]
\label{prop:fully-faithful-dynamics}
Let $\mathsf{CompHausDyn}$ be the category whose objects are pairs $(T,f)$ with
$T$ compact Hausdorff and $f:T\to T$ continuous, and whose morphisms
$(T,f)\to(T',g)$ are continuous maps $u:T\to T'$ satisfying
$u\circ f=g\circ u$. The assignment
\[
  (T,f)\longmapsto (X_T,X_f)
\]
defines a fully faithful functor from $\mathsf{CompHausDyn}$ to condensed
dynamical systems. Consequently, two compact Hausdorff dynamical systems are
topologically conjugate if and only if their associated condensed dynamical
systems are isomorphic in the category of condensed dynamical systems.
\end{proposition}

\begin{proof}
For compact Hausdorff spaces, the functor $T\mapsto X_T$ is fully faithful:
morphisms $X_T\to X_{T'}$ in $\Cond$ are exactly continuous maps
$T\to T'$ \cite{ClausenScholze,LeStumCondensed}. Under this identification, the
condition that a morphism $X_T\to X_{T'}$ commute with $X_f$ and $X_g$ is exactly
the condition that the corresponding continuous map $u:T\to T'$ satisfy
$u\circ f=g\circ u$. Hence the induced functor on dynamical systems is fully
faithful. Isomorphisms in this category are precisely topological conjugacies.
\end{proof}

\section{Fixed, periodic and preperiodic loci}
\label{sec:periodic}

We first recall the elementary construction behind the notation. If
$a,b:A\to B$ are two maps of sets, their \emph{equalizer} is
\[
  \Eq(a,b)=\{x\in A:a(x)=b(x)\}.
\]
It is the part of $A$ on which the two processes $a$ and $b$ give the same
answer. Categorically it is displayed as
\[
\begin{tikzcd}
\Eq(a,b) \arrow[r, hook] & A \arrow[r, shift left=0.45ex, "a"]
  \arrow[r, shift right=0.45ex, swap, "b"] & B
\end{tikzcd}
\]
The same construction exists in any topos. In dynamics, the two maps are
\[
  a=\mathrm{id}_X,\qquad b=f^n,
\]
so the equalizer is precisely the locus where $f^n$ agrees with the identity.
The ordinary set of such points is recovered by evaluation at the terminal
profinite set. The internal construction is made first as a subobject.

We also need the subobject notation in $\Cond$. Let $X$ be a condensed set. A
subobject of $X$ is an isomorphism class of
monomorphisms
\[
  A\hookrightarrow X
\]
in $\Cond$. We write $\Sub_{\Cond}(X)$ for the partially ordered set of
subobjects of $X$, ordered by inclusion: $A\le B$ if the inclusion
$A\hookrightarrow X$ factors through $B\hookrightarrow X$.

Because $\Cond$ is a Grothendieck topos, the subobject poset
$\Sub_{\Cond}(X)$ has the lattice operations needed here
\cite{MacLaneMoerdijk,JohnstoneTopos}. It is a complete lattice, so arbitrary
meets $\bigwedge_i A_i$ and joins $\bigvee_i A_i$ exist. The meet is
intersection, computed by pullback. The join is the least subobject containing
the given family. If $u:Y\to X$ is a morphism
in $\Cond$, inverse image of subobjects is defined by pullback and preserves
finite meets and arbitrary joins. Appendix~\ref{app:categorical-background}
gives the diagrams for these operations. Thus the formulas below are literal
operations in the subobject lattice of the phase object $X$.

\begin{definition}
Let $(X,f)$ be a condensed dynamical system.
\begin{enumerate}[label=(\alph*)]
\item The \emph{fixed locus} is the subobject
\[
  \Fix(X,f)=\Eq(\mathrm{id}_X,f)\hookrightarrow X.
\]
\item For $n\ge1$, the \emph{$n$-periodic locus} is the subobject
\[
  \Per_n(X,f)=\Eq(\mathrm{id}_X,f^n)\hookrightarrow X.
\]
It is the locus of points whose period divides $n$.
\item The \emph{periodic locus} is the join
\[
  \Per(X,f)=\bigvee_{n\ge1}\Per_n(X,f)\hookrightarrow X.
\]
\item The \emph{preperiodic locus} is
\[
  \PrePer(X,f)=\bigvee_{r\ge0}f^{-r}\Per(X,f)\hookrightarrow X,
\]
where inverse images of subobjects are formed by pullback.
\end{enumerate}
\end{definition}

\begin{figure}[htbp]
\centering
\begin{tikzcd}
\Per_n(X,f) \arrow[r, hook] &
X \arrow[r, shift left=0.45ex, "\mathrm{id}_X"]
  \arrow[r, shift right=0.45ex, swap, "f^n"] &
X
\end{tikzcd}
\[
  \Per_n(X,f)(S)
  =
  \{x\in X(S):f_S^n(x)=x\}.
\]
\caption{The $n$-periodic locus as an equalizer.}
\label{fig:periodic-equalizer}
\end{figure}

The divisibility poset on positive integers is directed, and the inclusions
\[
  \Per_m(X,f)\hookrightarrow \Per_n(X,f)
  \qquad(m\mid n)
\]
come from the implication $f^m(x)=x\Rightarrow f^n(x)=x$. Thus the periodic
locus may equivalently be written as the filtered colimit over the divisibility
poset.

Classical cycle counts start from the same fixed loci of iterates. For a finite
dynamical system, set
\[
  a_n=|\Fix(f^n)|=|\Per_n(f)|.
\]
If $b_n$ denotes the number of points of exact period $n$, then
\[
  a_n=\sum_{d\mid n} b_d,
  \qquad
  b_n=\sum_{d\mid n}\mu(d)a_{n/d},
\]
where $\mu$ is the M\"obius function, and the number of cycles of exact length
$n$ is $c_n=b_n/n$. This is the usual
cycle-counting viewpoint in finite and $p$-adic dynamics
\cite{KhrennikovNilsonPadicDynamics,SilvermanADS}. Here the same fixed loci
$\Per_n(X,f)$ are retained as functorial objects before the passage to numbers.
The primitive operations used here are joins, meets and inverse images.
Classical exact-period strata are obtained by removing lower-period loci.
Internally, this requires the relevant lower-period subobject to be
complemented, which need not hold in a general topos. The period-dividing
equalizers $\Per_n(X,f)$ are therefore the primary internal periodic objects.

\begin{remark}[Local boundedness of periods]
A section of $\Per(X,f)$ over a profinite set $S$ is locally represented,
after passing to a finite jointly surjective cover, by sections lying in
finitely many $\Per_{n_i}(X,f)$. Replacing the $n_i$ by their least common
multiple gives one period on that cover. Every $S$-section therefore has a
common period after a finite cover, without a single period fixed uniformly for
all profinite sets and all sections. This is exactly the finite-stage detection of
Lemma~\ref{lem:finite-stage-detection} applied to the filtered union of the
$\Per_n(X,f)$ over the divisibility poset.
\end{remark}

\begin{proposition}[Classical recovery]
\label{prop:classical-recovery}
Let $(T,f)$ be a topological dynamical system with $T$ compactly generated weak
Hausdorff, and let $(X_T,X_f)$ be its associated condensed dynamical system.
Let $a:*\to X_T$ correspond to a point $a\in T$.
Then
\begin{align*}
  \Fix(X_T,X_f)(*)
  &\cong
  \{x\in T:f(x)=x\},
\\
  \Per_n(X_T,X_f)(*)
  &\cong
  \{x\in T:f^n(x)=x\},
\\
  \Hist^{-}(X_T,X_f)(*)
  &\cong
  \{(x_i)_{i\ge0}:f(x_{i+1})=x_i\ \text{for all }i\ge0\},
\\
  T_n(a)(*)
  &\cong
  \{x\in T:f^n(x)=a\}.
\end{align*}
Similarly, $\Per(X_T,X_f)(*)$ and $\PrePer(X_T,X_f)(*)$ recover the usual
periodic and preperiodic points of $f$.
\end{proposition}

\begin{proof}
Since $X_T(*)=\Cont(*,T)\cong T$, evaluation of the equalizers
$\Eq(\mathrm{id}_{X_T},X_f)$ and $\Eq(\mathrm{id}_{X_T},X_f^n)$ at $*$ gives
the ordinary fixed and $n$-periodic loci. For $\Per$, Lemma~\ref{lem:finite-stage-detection}
shows that a terminal section of the filtered join is represented globally in
some finite stage $\Per_N(X_T,X_f)(*)$. Thus such a section is a classically
periodic point, and a classical periodic point with $f^n(x)=x$ defines a section
of $\Per_n(X_T,X_f)(*)$, hence of $\Per(X_T,X_f)(*)$. Since the image under $f$
of a periodic point is periodic, the subobjects
$f^{-r}\Per(X,f)$ form an increasing family. Applying the same finite-stage
lemma to
$\PrePer(X,f)=\bigvee_{r\ge0}f^{-r}\Per(X,f)$ shows that every terminal section
lies in $f^{-r}\Per(X,f)(*)$ for some $r$. These are precisely the classical
preperiodic points.
The formula for $\Hist^{-}(X_T,X_f)(*)$ follows because limits of sheaves are
computed sectionwise. The formula for $T_n(a)(*)$ is the ordinary pullback
description of the fiber of $f^n$ over $a$.
\end{proof}

\begin{proposition}[Invariance]
\label{prop:invariance}
For any condensed dynamical system $(X,f)$, each $\Per_n(X,f)$ is forward
invariant. Hence $\Per(X,f)$ is forward invariant. The preperiodic locus
$\PrePer(X,f)$ is both forward invariant and backward invariant in the sense of
Definition~\ref{def:invariance}. In general $\Per_n(X,f)$ need not be backward
invariant. The example after the proof shows this.
\end{proposition}

\begin{proof}
If $x$ is a section of $\Per_n(X,f)$, then $f^n(x)=x$, so
\[
  f^n(f(x))=f(f^n(x))=f(x).
\]
Thus $f(x)$ is again a section of $\Per_n(X,f)$. Since $\Per(X,f)$ is the join
of the $\Per_n(X,f)$, it is forward invariant.

For preperiodicity, if $f^r(x)$ is periodic, then $f^r(f(x))$ is periodic, so
the locus is forward invariant. Backward invariance follows formally because
inverse image preserves joins:
\begin{align*}
  f^{-1}\PrePer(X,f)
  &=
  f^{-1}\left(\bigvee_{r\ge0}f^{-r}\Per(X,f)\right)
\\
  &=
  \bigvee_{r\ge0}f^{-(r+1)}\Per(X,f)
  \le
  \PrePer(X,f).
\end{align*}
\end{proof}

For the last assertion, take the map $\{0,1\}\to\{0,1\}$ with $0\mapsto1$ and
$1\mapsto1$. Its fixed locus is $\{1\}$, while the inverse image of that fixed
locus is all of $\{0,1\}$.

\begin{proposition}[Limit compatibility]
\label{prop:limit}
Let $(X_i,f_i)$ be a projective system of condensed dynamical systems with
limit $(X,f)$. Then, for every $n\ge1$, there is a canonical isomorphism
\[
  \Per_n(X,f)\cong\varprojlim_i\Per_n(X_i,f_i).
\]
\end{proposition}

\begin{proof}
The limit system is obtained by taking $X=\varprojlim_i X_i$ and
$f=\varprojlim_i f_i$. Hence $f^n=\varprojlim_i f_i^n$, because composition is
compatible with the transition maps of the projective system. The object
$\Per_n(X,f)$ is the equalizer of $\mathrm{id}_X$ and $f^n$. Since equalizers are
limits and limits commute with limits in any complete category, this equalizer is
the projective limit of the equalizers
$\Eq(\mathrm{id}_{X_i},f_i^n)=\Per_n(X_i,f_i)$.
\end{proof}

\begin{remark}
The proposition is categorical. It says that fully compatible periodic points
in an inverse limit are exactly compatible systems of finite-level periodic
points. Lifting a periodic point visible at one fixed finite level to the
inverse limit is the separate arithmetic question addressed, in favorable
cases, by Hensel-type arguments.
\end{remark}

\begin{example}[Periodic families]
\label{ex:periodic-families}
Let $P_n=\{x\in T:f^n(x)=x\}$ be a finite discrete periodic locus for a
classical system $(T,f)$, and suppose the condensed periodic locus is
represented by $P_n$. If $S$ is finite discrete, then
\[
  \Per_n(X_T,X_f)(S)\cong\operatorname{Map}(S,P_n).
\]
A section is an $S$-tuple of $n$-periodic points.

If $S$ is profinite, then
\[
  \Per_n(X_T,X_f)(S)\cong\Cont(S,P_n),
\]
the set of locally constant maps from $S$ to $P_n$. Equivalently, such a
section is a labeling of a finite clopen partition of $S$ by $n$-periodic
points. The number $|P_n|$ is still visible, but the condensed object also
remembers how periodic data pull back along maps of profinite sets.
\end{example}

\subsection{Arithmetic towers and Hensel lifting}
\label{sec:arithmetic}

The preceding constructions become concrete in a congruence tower. A polynomial
map on $\Zp$ has finite reductions modulo $p^m$, and these finite systems are
compatible as $m$ varies. The condensed periodic locus packages that tower into
one object. Hensel's lemma then answers a different question: when does a
particular solution seen modulo $p$ lift to an actual $p$-adic periodic point?

Fix a prime $p$ and a polynomial $f\in\Zp[x]$. The reductions modulo $p^m$
define compatible finite systems
\[
  f_m:\Z/p^m\Z\to\Z/p^m\Z.
\]
Their inverse limit is the continuous system
\[
  f:\Zp\to\Zp.
\]
The corresponding condensed phase object is the inverse limit of the finite
condensed phase objects:
\[
  X_{\Zp}\cong \varprojlim_m X_{\Z/p^m\Z}.
\]
Indeed, after evaluation on a profinite set \(S\), this is the usual identity
\[
  \Cont\!\left(S,\varprojlim_m\Z/p^m\Z\right)
  \cong
  \varprojlim_m\Cont(S,\Z/p^m\Z),
\]
valid because continuous maps into a limit are the same as compatible maps into
its finite stages.
Applying Proposition~\ref{prop:limit} gives
\[
  \Per_n(\Zp,f)\cong
  \varprojlim_m \Per_n(\Z/p^m\Z,f_m),
\]
where $\Zp$ and $\Z/p^m\Z$ are viewed as condensed sets attached to their
natural profinite topologies.

The compatibility of periodic and preimage data across finite levels is a
recurring theme in arithmetic dynamics, for instance in dynamical
N\'eron--Ogg--Shafarevich criteria formulated through orbital arboreal
representations. See \cite{PerezBuendiaNOS} for one such arithmetic context.
Here the elementary Hensel-lifting mechanism supplies the local arithmetic
input.

\begin{example}[$x\mapsto x^2$ on $\Zp$]
Let $p$ be odd and let $f(x)=x^2$. The equation $f^n(x)=x$ is
\[
  x^{2^n}=x,
  \qquad\text{equivalently}\qquad
  x(x^{2^n-1}-1)=0.
\]
In $\Zp$, the solutions are
\[
  \Per_n(\Zp,f)=\{0\}\cup\mu_d,
  \qquad d=\gcd(2^n-1,p-1),
\]
where $\mu_d\subset\Zp^\times$ denotes the $d$-th roots of unity, represented
by Teichmüller lifts. This uses the standard description of the prime-to-$p$
torsion in $\Zp^\times$ for odd $p$ \cite{GouveaPadic}. At finite congruence
levels, non-liftable solutions may appear for those $n$ with
$p\mid(2^n-1)$. Indeed, at a nonzero root the derivative of
$x^{2^n-1}-1$ is $(2^n-1)x^{2^n-2}$, so Hensel's nonvanishing condition may
fail modulo $p$. The inverse limit retains only the compatible $p$-adic
solutions.

For instance, take $p=7$ and $n=3$. Then the equation is $x(x^7-1)=0$. A direct
count gives
\[
\begin{array}{c|c|c}
m & |\Per_3(\Z/7^m\Z,x^2)| &
\operatorname{im}\bigl(\Per_3(\Z/7^m\Z,x^2)\to
\Per_3(\Z/7^{m-1}\Z,x^2)\bigr)
\\ \hline
1 & 2 & \text{--}\\
2 & 8 & \{0,1\}\subset \Z/7\Z\\
3 & 8 & \{0,1\}\subset \Z/49\Z.
\end{array}
\]
At level \(49\), for example, the solutions are
\(\{0,1,8,15,22,29,36,43\}\), and all the nonzero ones reduce to \(1\) modulo
\(7\). At level \(343\), all eight solutions reduce modulo \(49\) to \(0\) or
\(1\). Thus many finite-level solutions do not form new compatible
\(7\)-adic branches. The inverse limit keeps only the two compatible solutions
\(\{0,1\}\).
\end{example}

\begin{proposition}[Hensel lifting for periodic points]
\label{prop:hensel}
Let $f\in\Zp[x]$, let $n\ge1$, and set
\[
  g(x)=f^n(x)-x.
\]
Suppose that $a_0\in\Z/p\Z$ satisfies
\[
  g(a_0)\equiv0\pmod p,
  \qquad
  g'(a_0)\not\equiv0\pmod p.
\]
Then there exists a unique $a\in\Zp$ with $a\equiv a_0\pmod p$ and
$f^n(a)=a$.
\end{proposition}

\begin{proof}
This is Hensel's lemma applied to the polynomial $g(x)=f^n(x)-x$. See
\cite{GouveaPadic} for the basic $p$-adic statement and
\cite[Ch.~2]{SilvermanADS} for the arithmetic-dynamical context.
\end{proof}

\begin{remark}
Condensed mathematics organizes the compatible tower of periodic loci.
Hensel's lemma decides when a specific finite-level solution lifts uniquely.
These are complementary mechanisms.
\end{remark}

\section{Discrete and compact periodic structures}
\label{sec:not-counting}

The internal object $\Per_n(X,f)$ contains more information than the number of
ordinary periodic points or cycles. The ordinary count of points fixed by
$f^n$ is $|\Per_n(X,f)(*)|$, obtained by evaluating the $n$-periodic locus at
the terminal profinite set. The sheaf
\[
  S\longmapsto \Per_n(X,f)(S)
\]
also records profinite families, functorial pullback in $S$, and finite-level
compatibility. In this sense, condensed dynamics refines the cycle-counting
viewpoint: it retains the family structure from which the usual numerical
counts are obtained.

If $P_n$ is a finite discrete periodic set, then for a finite discrete profinite set
$S$ one has
\[
  \Per_n(X,f)(S)\cong\mathrm{Map}(S,P_n),
  \qquad
  |\Per_n(X,f)(S)|=|P_n|^{|S|}.
\]
At this finite numerical level, the count is determined by $|P_n|$. The
additional structure appears in four related ways.

\begin{enumerate}[label=(\roman*)]
\item A morphism $S'\to S$ of profinite sets pulls an $S$-family back to an
$S'$-family. This functoriality is absent from cardinality alone.

\item Time parameters in the action morphism belong to $\Ndisc(S)$. They are
locally constant maps $S\to\N$. Compactness of $S$ makes each one constant on
the pieces of a finite clopen partition.

\item Orbit quotients are formed in the sheaf category. The object
$X_{\N}=\Coeq(\mathrm{id}_X,f)$ includes the sheafification required by the
topos. Its value on $S$ therefore retains descent data and need not equal the
pointwise quotient $X(S)/\langle x\sim f_S(x)\rangle$.

\item Infinite profinite parameter spaces detect topology. For a finite
discrete target $P_n$, the set $\Cont(S,P_n)$ consists of locally constant maps.
For a compact Cantor target $K$, the set $\Cont(S,K)$ also contains continuous
families that do not factor through a finite periodic stratum.
\end{enumerate}

\begin{definition}
\label{def:ind}
Let $P$ be a small set. We write $P_{\mathrm{ind}}$ for the filtered
colimit, in $\Cond$, of the finite subsets of $P$, each viewed as a finite
discrete condensed set, and call it the \emph{ind-discrete object} on $P$. Thus
the notation emphasizes its filtered presentation by finite discrete subsets.
\end{definition}

Lemma~\ref{lem:finite-stage-detection} describes how sections of this filtered
presentation are controlled by finite stages. For a set $P$, let
$P_{\mathrm{disc}}$ denote $P$ with the discrete topology. A \emph{Cantor
space} is a nonempty compact metrizable, perfect and totally disconnected
space. Every such space is homeomorphic to $\{0,1\}^{\N}$.

\begin{proposition}[Discrete realization and dense Cantor separation]
\label{prop:discrete-compact}
\leavevmode
\begin{enumerate}[label=\textnormal{(\alph*)},leftmargin=2.2em]
\item The canonical morphism $P_{\mathrm{ind}}\to X_{P_{\mathrm{disc}}}$ is an
isomorphism for every small set $P$.
\item If $K$ is a Cantor space and $P\subset K$ is countably infinite and dense,
then the natural morphism $P_{\mathrm{ind}}\to X_K$ is not an isomorphism in
$\Cond$.
\end{enumerate}
\end{proposition}

\begin{proof}
Let $S$ be profinite. Every continuous map $S\to P_{\mathrm{disc}}$ has finite
image because $S$ is compact and $P_{\mathrm{disc}}$ is discrete. It therefore
factors through a finite subset $F\subset P$. Hence
\[
  \varinjlim_{F\subset P,\,F\text{ finite}}\Cont(S,F)
  \cong \Cont(S,P_{\mathrm{disc}}).
\]
The presheaf colimit on the left is therefore $X_{P_{\mathrm{disc}}}$ and is
already a sheaf. Its sheafification is unchanged, which proves (a).

For (b), let $K$ be a Cantor space and let $P\subset K$ be
countably infinite and dense.
Evaluate at a Cantor profinite set $S$, for instance $S=\{0,1\}^{\N}$.
By (a), $P_{\mathrm{ind}}(S)=\Cont(S,P_{\mathrm{disc}})$ is the set of locally
constant maps $S\to P$. The space $S$ has a countable clopen basis. Every
clopen subset is compact and hence a finite union of basis elements. It follows
that $S$ has only countably many clopen subsets and only countably many finite
clopen partitions. Each locally constant map is determined by one such
partition and a finite list of values in $P$. Since $P$ is countable,
$\Cont(S,P_{\mathrm{disc}})$ is countable. The set
$X_K(S)=\Cont(S,K)$ is uncountable because it contains one constant map for
each point of $K$. It also contains sections that are not locally constant,
for instance a homeomorphism $S\cong K$. Thus evaluation at
$S$ distinguishes the two objects, and the natural morphism
$P_{\mathrm{ind}}\to X_K$ cannot be an isomorphism in $\Cond$.
\end{proof}

\subsection{A Cantor-type Julia application}
\label{sec:julia-cantor}

We now apply the preceding separation to a classical non-Archimedean Julia set,
regarded first as a compact topological space and then as a condensed set. We
work with quadratic bad-reduction examples that admit symbolic models; see
\cite[Example~5.29 and Exercise~8.19]{BenedettoNAM}. Background on Berkovich
dynamics and repelling periodic points can be found in
\cite[Ch.~12]{BenedettoNAM} and \cite{BakerRumely}. The general density question
for repelling periodic points in classical non-Archimedean Julia sets remains
open in this form \cite{OkuyamaRepelling}.

In this section $J(f)$ denotes the classical Julia set of a rational map over
$\Cp$, equipped with the topology supplied by the usual non-Archimedean theory.
Its associated condensed set is
\[
  X_{J(f)}(S)=\Cont(S,J(f)).
\]
The observed structure described in the Introduction provides a setting for
defining Fatou and Julia loci by equicontinuity inside a topos. Here the
classical Julia set is used as an external topological input.

One concrete family is $f(z)=z^2+c$ over $\Cp$, where $p$ is odd and
$|c|_p>1$. Its critical orbit escapes to infinity. If $\alpha$ is a fixed point,
translation by $\alpha$ conjugates $f$ to $z^2+\lambda z$, where
$\lambda=2\alpha$ and $|\lambda|_p>1$. The inverse-disk construction and
symbolic coding in \cite[Example~5.29 and Exercise~8.19]{BenedettoNAM} give a
topological conjugacy between the restriction of $f$ to its classical Julia
set and the one-sided full shift on two symbols. In particular, $J(f)$ is a
compact totally disconnected perfect set, hence a Cantor-type profinite space.
The metric picture is the following. Over $\Cp$ the absolute value is ultrametric, so
two balls are either nested or disjoint, never partially overlapping
(Appendix~\ref{app:topological-background}). On the chaotic locus the map then
splits the space into disjoint pieces at every scale, the branching pattern of a
Cantor set. Periodic sequences are dense in the shift space. Hence periodic
points are dense in $J(f)$. Write
$\Per_J(f)$ for the set of periodic points lying in the Julia set. Then
\[
  \overline{\Per_J(f)}=J(f),
\]
where the closure is taken in the classical topology.
Since the shift has only finitely many points fixed by each iterate, the set
\(\Per_J(f)\) is countable.

Let \(X_{f|J}:X_{J(f)}\to X_{J(f)}\) be the condensed endomorphism induced by
the restriction of \(f\) to \(J(f)\). For each \(n\ge1\),
\[
  \Per_n(X_{J(f)},X_{f|J})
  \cong
  X_{\Fix((f|_{J(f)})^n)}.
\]
Indeed, after evaluation on a profinite set \(S\), both sides consist of the
continuous maps \(\alpha:S\to J(f)\) satisfying
\((f|_{J(f)})^n\circ\alpha=\alpha\). The symbolic model makes each
\(\Fix((f|_{J(f)})^n)\) finite. Conversely, every finite subset of
\(\Per_J(f)\) is contained in \(\Fix((f|_{J(f)})^N)\) for some \(N\), by taking
a common multiple of the periods of its points. Hence
\[
  \Per(X_{J(f)},X_{f|J})
  \cong
  (\Per_J(f))_{\mathrm{ind}}
  \cong
  X_{(\Per_J(f))_{\mathrm{disc}}},
\]
where the second isomorphism is Proposition~\ref{prop:discrete-compact}(a).

We write
\[
  \Per_{J,\mathrm{ind}}(f):=(\Per_J(f))_{\mathrm{ind}}
\]
for this internal periodic locus, in the notation of Definition~\ref{def:ind}.
The compatible inclusions of finite subsets \(F\subset \Per_J(f)\) into \(J(f)\)
induce a natural monomorphism
\[
  \Per_{J,\mathrm{ind}}(f)\hookrightarrow X_{J(f)}.
\]
The morphism is a monomorphism because it evaluates on every profinite set as
the injective map induced by $\Per_J(f)\hookrightarrow J(f)$. At the terminal
profinite set it recovers the classically dense inclusion
$\Per_J(f)\subset J(f)$. It is nevertheless not an isomorphism in $\Cond$:
\[
  \Per_{J,\mathrm{ind}}(f)
  \not\simeq
  X_{J(f)}.
\]
The hypotheses of Proposition~\ref{prop:discrete-compact}(b) hold because
$\Per_J(f)$ is countable and dense and $J(f)$ is a Cantor space.

Thus the abstract separation becomes a familiar non-Archimedean picture:
periodic points are dense in the Julia set, but the dense periodic object is
still not the compact condensed object that contains its closure.

\subsection{Further directions}
\label{sec:outlook}

The constructions above provide the internal base layer for several adjacent
developments.

\begin{enumerate}[label=(\roman*)]
\item Numerical invariants, including fixed-point and cycle counts, may be
attached to finite-level evaluations of $\Per_n(X,f)$. Arithmetic refinements
along towers such as $\Z/p^m\Z$ belong to this numerical layer.

\item The observed Fatou--Julia theory developed in
\cite{PerezBuendiaFatuTopos} equips $(X,f)$ with an observable frame and an
admissible coverage. These data support cover-wise equicontinuity and observable
Julia loci.

\item Linear variants place related invariants in condensed abelian groups and
other linear condensed categories.

\item Replacing the discrete monoid $\Ndisc$ by more general monoid objects
in $\Cond$ leads to a separate theory of internal monoid actions, stabilizers and
quotients.
\end{enumerate}

Each direction begins with the orbit, history, preimage and periodic objects
constructed above. Global sections recover ordinary dynamical data, while
evaluation on profinite spaces preserves their variation across finite levels.
The condensed formalism holds these returns together as \emph{memorias},
gathered by the rising sea.

\appendix

\section{Categorical, site and topos background}
\label{app:categorical-background}

This appendix fixes the categorical notation and complements the standard
references on sites, topoi, Stone spaces and condensed sets
\cite{MacLaneMoerdijk,JohnstoneTopos,JohnstoneStoneSpaces,ClausenScholze,LeStumCondensed}.
It gives the definitions and permanence properties used in
Sections~\ref{sec:condensed} and~\ref{sec:periodic}.

\subsection{Sites and Grothendieck topoi}

\begin{definition}[Grothendieck pretopology and site]
Let $\mathcal C$ be a category with pullbacks. A \emph{Grothendieck
pretopology} on $\mathcal C$ assigns to each object $U$ a class of families
\[
  \{U_i\to U\}_{i\in I}
\]
called \emph{covering families}, subject to the following three axioms.
\begin{enumerate}[label=(\roman*)]
\item If $V\to U$ is an isomorphism, then $\{V\to U\}$ is a covering family.
\item If $\{U_i\to U\}_{i\in I}$ is a covering family and $V\to U$ is any
arrow, then the base-changed family
\[
  \{V\times_U U_i\to V\}_{i\in I}
\]
is a covering family.
\item If $\{U_i\to U\}_{i\in I}$ is a covering family and, for each $i$,
$\{V_{ij}\to U_i\}_{j\in J_i}$ is a covering family, then the composite family
\[
  \{V_{ij}\to U_i\to U\}_{i,j}
\]
is a covering family.
\end{enumerate}
A category equipped with such a pretopology is a \emph{site}. We use this
pretopological presentation throughout. The associated Grothendieck topology is
the one generated by these covers. For the standard formulation, see
\cite[Ch.~III]{MacLaneMoerdijk} and \cite[Vol.~2, C2]{JohnstoneTopos}.
\end{definition}

\begin{definition}[Sheaf condition]
Let $(\mathcal C,\tau)$ be a site. A presheaf is a contravariant functor
\[
  F:\mathcal C^{\mathrm{op}}\to\Set.
\]
It is a \emph{sheaf} if, for every covering family $\{U_i\to U\}_{i\in I}$, the
canonical diagram
\[
  F(U)\longrightarrow \prod_i F(U_i)
  \rightrightarrows
  \prod_{i,j}F(U_i\times_U U_j)
\]
is an equalizer. Equivalently, a section over $U$ is the same as a family of
sections over the $U_i$ whose two restrictions to each overlap
$U_i\times_U U_j$ agree. The category of sheaves on the site is denoted
$\Sh(\mathcal C,\tau)$, or simply $\Sh(\mathcal C)$ when
$\tau$ is fixed.
\end{definition}

\begin{definition}[Grothendieck topos]
A \emph{Grothendieck topos} is a category equivalent to
$\Sh(\mathcal C,\tau)$ for some site $(\mathcal C,\tau)$. The categorical
properties used in the main text are that a Grothendieck topos has all small
limits and colimits, limits of sheaves are computed sectionwise, subobject
posets are complete lattices, and colimits of sheaves are obtained from
presheaf colimits by sheafification. See
\cite{MacLaneMoerdijk,JohnstoneTopos}.
\end{definition}

\begin{definition}[The profinite site for condensed sets]
Fix a set-theoretic universe, or equivalently a small skeleton of the profinite
site, and let $\ProFin$ be the corresponding category of profinite spaces and
continuous maps. The finite jointly surjective pretopology $\tau_{\mathrm{fs}}$
on $\ProFin$ declares a finite family
\[
  \{U_i\to S\}_{i=1}^r
\]
to be a cover when the induced map of profinite spaces
\[
  \coprod_{i=1}^r U_i\longrightarrow S
\]
is surjective. A \emph{condensed set}, in the presentation used here, is a
sheaf
\[
\begin{tikzcd}[column sep=large]
\ProFin^{\mathrm{op}}
\arrow[r, "{X}"]
&
\Set ,
\end{tikzcd}
\]
for this pretopology. Thus
\[
  \Cond=\Sh(\ProFin,\tau_{\mathrm{fs}}).
\]
This is a Grothendieck topos \cite{ClausenScholze,LeStumCondensed}.
\end{definition}

For a condensed set $X$, the sheaf condition says that every finite jointly
surjective cover $\{U_i\to S\}$ gives an equalizer
\[
  X(S)\longrightarrow \prod_i X(U_i)
  \rightrightarrows
  \prod_{i,j}X(U_i\times_S U_j).
\]
The fiber product $U_i\times_S U_j$ is the overlap of the two profinite pieces
over $S$.

We use two standard facts about sheaf categories. Limits of sheaves are
computed as the corresponding limits of presheaves, hence sectionwise.
Colimits are obtained by taking the presheaf colimit and then sheafifying. This
is why inverse-limit constructions such as backward histories evaluate
sectionwise, while quotient constructions such as orbit quotients should not
be identified blindly with pointwise quotients.

\subsection{Subobjects, meets and joins}

Let $\mathcal E$ be a Grothendieck topos and let $X\in\mathcal E$. A subobject
of $X$ is an isomorphism class of monomorphisms $A\hookrightarrow X$. These
subobjects form a partially ordered set
\[
  \Sub_{\mathcal E}(X).
\]
The order is factorization: $A\le B$ if $A\hookrightarrow X$ factors through
$B\hookrightarrow X$.

The poset $\Sub_{\mathcal E}(X)$ is a complete lattice
\cite{MacLaneMoerdijk,JohnstoneTopos}. We use the following operations:
\[
\begin{array}{ccl}
A\wedge B & = & \text{intersection, computed by pullback},\\[2mm]
\displaystyle\bigvee_i A_i & = & \text{least subobject containing all }A_i,\\[2mm]
0_X & = & \text{initial subobject of }X,\\[1mm]
1_X & = & X\hookrightarrow X.
\end{array}
\]
The meet diagram is
\[
\begin{tikzcd}
A\wedge B \arrow[r] \arrow[d] & B \arrow[d, hook]\\
A \arrow[r, hook] & X
\end{tikzcd}
\]
The join can be described as the image of the canonical map
\[
\begin{tikzcd}
\coprod_i A_i \arrow[r] & X
\end{tikzcd}
\]
when the coproduct is formed in $\mathcal E$.

\begin{remark}
The subobject lattice is in fact a complete Heyting algebra. Its implication
$A\Rightarrow B$ is the largest subobject $C\hookrightarrow X$ such that
$C\wedge A\le B$. The main text uses meets, joins and inverse images, so no
Heyting implication is needed there.
\end{remark}

If $u:Y\to X$ is a morphism, inverse image of subobjects is defined by
pullback:
\[
\begin{tikzcd}
u^{-1}A \arrow[r] \arrow[d, hook] & A \arrow[d, hook]\\
Y \arrow[r, "u"] & X
\end{tikzcd}
\]
In a Grothendieck topos this inverse-image operation preserves finite meets
and arbitrary joins. This is the formal reason why expressions such as
\[
  \Per(X,f)=\bigvee_{n\ge1}\Per_n(X,f),
  \qquad
  \PrePer(X,f)=\bigvee_{r\ge0}f^{-r}\Per(X,f)
\]
are well-defined inside $\Cond$.

\subsection{The dynamical translation}

The following table gives the dynamical reading of the categorical operations.
\[
\begin{array}{c|c}
\text{Categorical operation} & \text{Dynamical meaning}\\
\hline
\Eq(\mathrm{id}_X,f) & \text{fixed points or fixed families}\\
\Eq(\mathrm{id}_X,f^n) & \text{points or families fixed by }f^n\\
A\wedge B & \text{satisfy both conditions}\\
\bigvee_i A_i & \text{satisfy one of the conditions locally}\\
f^{-1}(A) & \text{points whose image lies in }A\\
X(S) & S\text{-parameterized sections of }X
\end{array}
\]
These operations apply the internal logic of the ambient topos directly to the
endomorphism $f:X\to X$.

\section{Profinite and p-adic background}
\label{app:topological-background}

This appendix collects the elementary topological facts about profinite spaces,
$\Zp$, $\Qp$ and $\Cp$ used in the condensed constructions and examples.

\begin{definition}[Profinite space]
A \emph{profinite space} is a topological space isomorphic to an inverse limit
of finite discrete spaces:
\[
  S\simeq \varprojlim_i S_i,\qquad S_i\ \text{finite discrete}.
\]
Equivalently, it is a compact Hausdorff totally disconnected space. Its
topology has a basis of clopen subsets, and every continuous map from a
profinite space to a finite discrete set factors through a finite quotient.
For the classical theory of Stone and profinite spaces see
\cite{JohnstoneStoneSpaces}. For an arithmetic introduction, see the published
account \cite{PerezBuendiaStoneGaceta} and the longer self-contained version
\cite{PerezBuendiaStoneArxiv}.
\end{definition}

Profinite spaces serve as compact parameter spaces whose maps to finite targets
are controlled at finite levels. Evaluating a condensed set on
$S\in\ProFin$ keeps this compatibility data together with its terminal sections
$X(*)$.

\begin{definition}[$p$-adic integers]
The ring of $p$-adic integers is the inverse limit of the finite quotient
rings
\[
  \Zp=\varprojlim_m \Z/p^m\Z.
\]
It is a compact profinite ring. Its basic clopen balls are the residue classes
\[
  a+p^m\Zp=\{x\in\Zp:x\equiv a\pmod {p^m}\}.
\]
For fixed $m$, these balls form a finite partition of $\Zp$ indexed by
$\Z/p^m\Z$. Refining from level $m$ to level $m+1$ splits each ball into
$p$ smaller clopen balls. This is the topological content behind the inverse
limit description.
\end{definition}

\begin{figure}[htbp]
\centering
\resizebox{\textwidth}{!}{%
\begin{tikzpicture}[
  every node/.style={font=\small},
  ball/.style={draw=blue!55!black, rounded corners=2pt, inner sep=2.6pt,
               fill=blue!5, font=\small},
  leaf/.style={draw=blue!55!black, rounded corners=2pt, inner sep=2.2pt,
               fill=blue!3, font=\footnotesize},
  e/.style={-{Latex[length=1.6mm]}, gray!65},
]
\def\xax{-1.3}
\draw[gray!55,-{Latex[length=2mm]}] (\xax,0.7)--(\xax,-4.7);
\node[gray!55!black,font=\scriptsize,rotate=90] at (\xax-0.5,-2) {radius of the balls};
\foreach \y/\lab in {0/{3^{0}=1},-2/{3^{-1}},-4/{3^{-2}}}{
  \draw[gray!55] (\xax-0.12,\y)--(\xax+0.12,\y);
  \node[left,font=\scriptsize,gray!45!black] at (\xax-0.18,\y) {$\lab$};
  \draw[gray!25,dashed] (\xax+0.12,\y)--(16.4,\y);
}
\node[ball] (r) at (8,0) {$\mathbb{Z}_3$};
\node[ball] (a0) at (2,-2)  {$0+3\mathbb{Z}_3$};
\node[ball] (a1) at (8,-2)  {$1+3\mathbb{Z}_3$};
\node[ball] (a2) at (14,-2) {$2+3\mathbb{Z}_3$};
\foreach \a in {a0,a1,a2}{\draw[e] (r)--(\a);}
\node[leaf] (b0) at (0,-4)   {$0+9\mathbb{Z}_3$};
\node[leaf] (b1) at (2,-4)   {$3+9\mathbb{Z}_3$};
\node[leaf] (b2) at (4,-4)   {$6+9\mathbb{Z}_3$};
\node[leaf] (c0) at (6,-4)   {$1+9\mathbb{Z}_3$};
\node[leaf] (c1) at (8,-4)   {$4+9\mathbb{Z}_3$};
\node[leaf] (c2) at (10,-4)  {$7+9\mathbb{Z}_3$};
\node[leaf] (d0) at (12,-4)  {$2+9\mathbb{Z}_3$};
\node[leaf] (d1) at (14,-4)  {$5+9\mathbb{Z}_3$};
\node[leaf] (d2) at (16,-4)  {$8+9\mathbb{Z}_3$};
\foreach \a/\b in {a0/b0,a0/b1,a0/b2, a1/c0,a1/c1,a1/c2, a2/d0,a2/d1,a2/d2}{
  \draw[e] (\a)--(\b);}
\node[gray!55!black,font=\scriptsize] at (8,-4.85)
  {$\cdots$ each ball splits again into $3$ balls of radius $3^{-3}$, and so on $\cdots$};
\end{tikzpicture}}

\medskip

\resizebox{\textwidth}{!}{%
\begin{tikzpicture}[
  every node/.style={font=\footnotesize},
  lab/.style={font=\scriptsize},
]
\def\W{12}
\def\h{0.9}
\foreach \rr/\lab in {0/{3^{0}=1},1/{3^{-1}},2/{3^{-2}},3/{3^{-3}}}{
  \node[lab,gray!45!black,left] at (-0.25,-\rr*\h-\h/2) {$\lab$};
}
\node[lab,gray!55!black,rotate=90] at (-1.45,-1.6) {radius};
\draw[fill=blue!6,draw=blue!55!black] (0,0) rectangle (\W,-\h);
\node at (\W/2,-\h/2) {$\mathbb{Z}_3$};
\foreach \i/\lab in {0/{0+3\mathbb{Z}_3},1/{1+3\mathbb{Z}_3},2/{2+3\mathbb{Z}_3}}{
  \draw[fill=blue!10,draw=blue!55!black] (\i*\W/3,-\h) rectangle (\i*\W/3+\W/3,-2*\h);
  \node at (\i*\W/3+\W/6,-1.5*\h) {$\lab$};
}
\foreach \i/\lab in {0/0,1/3,2/6,3/1,4/4,5/7,6/2,7/5,8/8}{
  \draw[fill=blue!14,draw=blue!55!black] (\i*\W/9,-2*\h) rectangle (\i*\W/9+\W/9,-3*\h);
  \node[font=\scriptsize] at (\i*\W/9+\W/18,-2.5*\h) {$\lab$};
}
\foreach \i in {0,...,26}{
  \draw[fill=blue!18,draw=blue!40!black,line width=0.2pt]
        (\i*\W/27,-3*\h) rectangle (\i*\W/27+\W/27,-4*\h);
}
\draw[decoration={brace,amplitude=4pt},decorate,gray!60]
   (\W+0.15,0) -- (\W+0.15,-3*\h)
   node[midway,right,lab,gray!50!black,align=left,xshift=2pt]
   {a point of $\mathbb{Z}_3$\\ $=$ a nested\\ chain of balls};
\end{tikzpicture}}
\caption{Clopen balls in $\mathbb{Z}_p$ (shown for $p=3$): a rooted tree and the
same nested partitions.}
\label{fig:zp-ball-tree}
\end{figure}
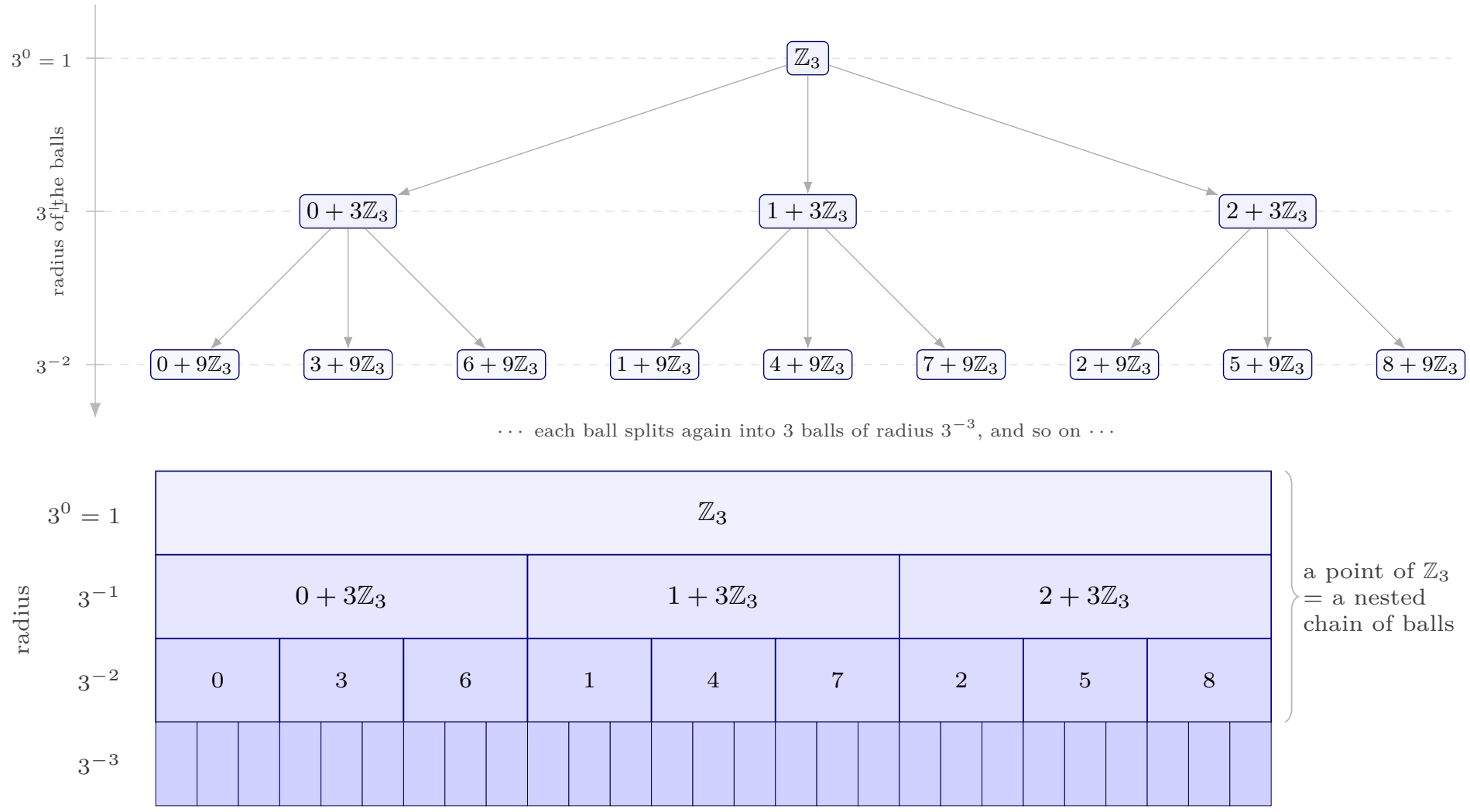

\begin{definition}[$p$-adic valuation, metric and fields]
For $x\in\Q^\times$, write $x=p^k(a/b)$ with $a,b\in\Z$ not divisible by $p$,
and set $v_p(x)=k$. Also set $v_p(0)=+\infty$. The $p$-adic absolute value on
$\Q$ is
\[
  |x|_p=p^{-v_p(x)}.
\]
It is non-Archimedean:
\[
  |x+y|_p\le \max\{|x|_p,|y|_p\}.
\]
The associated metric is $d_p(x,y)=|x-y|_p$. On $\Z$ this is
\[
  d_p(a,b)=
  \begin{cases}
  0, & a=b,\\
  p^{-v_p(a-b)}, & a\ne b.
  \end{cases}
\]
The completion of $\Q$ for $d_p$ is $\Qp$. Its valuation ring is $\Zp$, and
algebraically
\[
  \Qp\cong\operatorname{Frac}(\Zp).
\]
The field $\Cp$ is the completion of an algebraic closure of $\Qp$:
\[
  \Cp=\widehat{\Qpalg}.
\]
The following diagram relates the Archimedean and $p$-adic completions of
$\Q$, with $\Zp$ sitting as the valuation ring of $\Qp$:
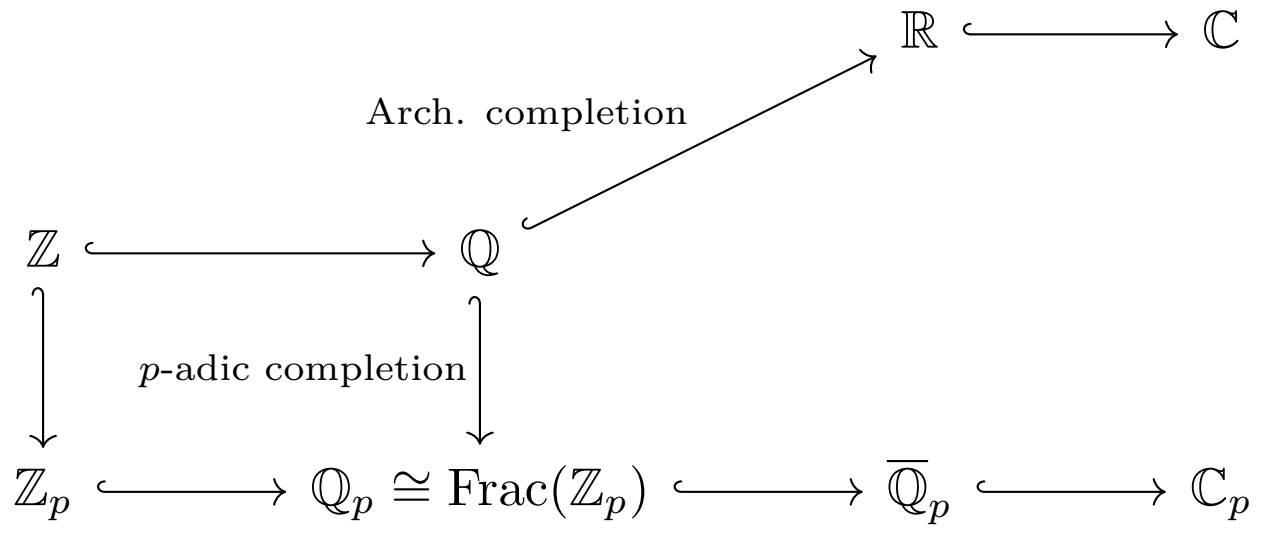
\[
\begin{tikzcd}[column sep=large,row sep=large]
& & \R \arrow[r, hook] & \C \\
\Z \arrow[r, hook] \arrow[d, hook'] &
\Q \arrow[ur, hook, "\text{Arch. completion}"] \arrow[d, hook',
  "\text{$p$-adic completion}"'] & & \\
\Zp \arrow[r, hook] &
\Qp\cong\operatorname{Frac}(\Zp) \arrow[r, hook] &
\Qpalg \arrow[r, hook] &
\Cp
\end{tikzcd}
\]
Here the maps $\Z\to\Zp$ and $\Q\to\Qp$ are the usual dense embeddings into the
$p$-adic completions, $\Zp$ is the valuation ring of $\Qp$, and the $p$-adic
absolute value extends uniquely along the lower branch.
For elementary background on $p$-adic absolute values, completions and the
topology of $\Qp$, see \cite{GouveaPadic}. For their use in arithmetic dynamics,
see \cite{SilvermanADS,BenedettoNAM}.
\end{definition}

The inverse-limit construction of $\Zp$ agrees with the analytic one. A Cauchy
sequence in $\Z$ for $d_p$ is eventually constant modulo each $p^m$, hence
determines a compatible system of residues. Conversely, such a compatible
system can be approximated by integers and defines a Cauchy sequence. This gives
the canonical isomorphism between the completion of $\Z$ for $d_p$ and
$\varprojlim_m\Z/p^m\Z$.

\begin{proposition}[Basic topology of $\Qp$]
\label{prop:basic-topology-q-p}
The field $\Qp$ is locally compact, Hausdorff and totally disconnected, but it
is not compact. It is obtained from $\Zp$ by inverting $p$, and topologically it
is the increasing union of compact open subsets
\[
  \Qp=\bigcup_{N\ge0}p^{-N}\Zp.
\]
Equivalently, $\Zp$ is the closed unit ball in $\Qp$, while $p^m\Zp$ is the
ball of radius $p^{-m}$ around $0$. Every ball in $\Qp$ is both open and
closed, and two balls are either disjoint or one contains the other. This
nested-or-disjoint behavior is a direct consequence of the ultrametric
inequality
\[
  |x+y|_p\le \max\{|x|_p,|y|_p\}.
\]
Indeed, the ultrametric inequality makes every point of a ball a center of it:
if $|x-a|_p\le r$, then any $y$ with $|y-x|_p\le r$ satisfies
$|y-a|_p\le\max\{|y-x|_p,|x-a|_p\}\le r$. Hence if two balls share a point, the
one of smaller radius lies inside the other. Every nonzero element $x\in\Qp$
can be written uniquely as
\[
  x=p^k u,\qquad k\in\Z,\quad u\in\Zp^\times,
\]
and then $v_p(x)=k$ and $|x|_p=p^{-v_p(x)}$. See \cite{GouveaPadic} for the
full development.
\end{proposition}

\[
  \Zp\subset p^{-1}\Zp\subset p^{-2}\Zp\subset\cdots.
\]

For condensed dynamics, compact totally disconnected spaces such as $\Zp$ and
Cantor-type Julia sets serve as parameter spaces for continuous families. The
nested chain of balls in Figure~\ref{fig:zp-ball-tree} is the metric counterpart
of compatibility along a congruence tower. Such a chain determines an
inverse-limit point. A general condensed section over $S$ may vary continuously
without factoring through a single finite stage.

\section*{Acknowledgments}

The author thanks the P-ADAGIO working group
(\url{https://www.cimat.mx/~rogelio.perez/fr/grupos-de-trabajo/p-adagio/}).
He also thanks the anonymous referee for a careful reading and for suggestions
that improved the exposition.

\end{document}